\documentclass[11pt]{amsart}
\usepackage{amsmath,amsthm,amssymb,amsfonts, bm}
\usepackage[applemac]{inputenc}
\usepackage[pdftex]{hyperref}
\numberwithin{equation}{section}
\makeindex
\newcommand{\supp}{\mathrm{supp}}

\newcommand{\p}{\mathfrak{p}}

\newcommand{\q}{\mathfrak{q}}
\newcommand{\C}{\mathbb{C}}
\newcommand{\m}{\mathfrak{m}}
\newcommand{\GL}{\mathrm{GL}}

\newcommand{\R}{\mathbb{R}}

\newcommand{\Z}{\mathbb{Z}}
\newcommand{\Q}{\mathbb{Q}}
\newcommand{\A}{\mathbb{A}}
\newcommand{\I}{\mathbb{I}}

\DeclareFontFamily{OT1}{rsfs}{}
\DeclareFontShape{OT1}{rsfs}{n}{it}{<-> rsfs10}{}
\DeclareMathAlphabet{\mathscr}{OT1}{rsfs}{n}{it}

\newtheorem{lemma}{Lemma}
\newtheorem{thm}{Theorem}
\newtheorem{prop}{Proposition}

\newtheorem{defn}{Definition}

\begin{document}

\author{Valentin Blomer}
\address{Universit\"at G\"ottingen, Mathematisches Institut,  Bunsenstr. 3-5, 37073 G\"ottingen} \email{blomer@uni-math.gwdg.de}

\author{Farrell Brumley}
\address{ Institut \'Elie Cartan, Universit\'e Henri Poincar\'e Nancy 1, BP 239, 54506 Vand{\oe}uvre Cedex, France} \email{farrell.brumley@iecn.u-nancy.fr}
\title[On the Ramanujan conjecture over number fields]{On the Ramanujan conjecture over number fields}

\thanks{The first author was supported by a Volks\-wagen Lichtenberg Fellowship and by a European Research Council (ERC) starting grant 258713.   The second author is supported by the ANR grant Modunombres.}


\begin{abstract} We extend to an arbitrary number field the best known bounds towards Ramanujan for the group $\GL_n$, $n=2,3,4$.  In particular, we present a technique which overcomes the analytic obstacles posed by the presence of an infinite group of units.
\end{abstract}


\maketitle

\section{Introduction}

\subsection{Statement of results}

Since Ramanujan \cite{Ra}, in 1916, stated his conjecture on the size of the coefficients $\tau(n)$ of the $\Delta(z)$ function, the task of bounding Fourier coefficients of modular forms has occupied a venerable position in analytic and algebraic number theory.  Deligne \cite{De}, famously, proved the Ramanujan conjecture for weight $k\geq 2$ holomorphic Hecke cusp forms, a result which was recently widely extended by Harris-Taylor \cite{HT}. That said, important natural generalizations of Ramanujan's conjecture remain open, the most notable of which is the Ramanujan-Petersson/Selberg conjecture on weight-zero Maa{\ss} forms for congruence subgroups $\Gamma$ of ${\rm SL}_2(\Bbb{Z})$.

The most general formulation of the Ramanujan conjecture is expressed through representation theory. Let $K$ be a number field with ring of adeles $\A$.  Let $\pi$ be a cuspidal automorphic representation of $\GL_n(\A)$ with unitary central character.  Fix an identification $\pi\simeq\otimes_v \pi_v$.  Then $\pi_v$ is an irreducible unitary generic representation of $\GL_n(K_v)$ and the Ramanujan conjecture is the assertion that $\pi_v$ is tempered.

A non-tempered representation $\pi_v$ can be described in the following way.  There exists a standard parabolic subgroup $P$ of $\GL_n(K_v)$ of type $(n_1,\ldots , n_r)$ with unipotent radical $U$, irreducible tempered representations $\tau_j$ of $\GL_{n_j}(K_v)$, and real numbers $\sigma_j$ satisfying $\sigma_1>\cdots >\sigma_r$ such that $\pi_v$ is equivalent to the fully induced representation ${\rm Ind}(\GL_n(K_v),P;{\boldsymbol\tau}[{\boldsymbol \sigma}])$.  Here ${\boldsymbol\tau}[{\boldsymbol \sigma}]$ is the representation of the group $M=P/U\simeq \GL_{n_1}\times\cdots\times \GL_{n_r}$ given by ${\boldsymbol\tau}[{\boldsymbol \sigma}]=\tau_r[\sigma_1]\otimes\cdots\otimes\tau_1[\sigma_r]$, and $\tau[\sigma]$ is the twisted representation $g\mapsto \tau(g)|\det g|_v^\sigma$.

The size of the parameters $\sigma_j$ allow one to quantitatively measure the failure of a given local representation to be tempered.  Note that since $\pi_v$ is unitary we have $\{\tau_j[\sigma_j]\}=\{\widetilde{\tau}_j[-\sigma_j]\}$ as sets, from which we deduce that $\max_j \sigma_j\leq\delta$ is equivalent to $\max_j|\sigma_j|\leq\delta$.

For any non-tempered $\pi_v$ appearing in $\pi\simeq\otimes_v\pi_v$ we restore the dependence of the parameters $\sigma_j$ on $\pi$ and the place $v$, writing $\sigma_\pi(v,1),\ldots ,\sigma_\pi(v,r)$.  For the rest of this paper we put $m(\pi,v)=\max_j|\sigma_{\pi}(v, j)|$ if $\pi_v$ is non-tempered, and $m(\pi,v)=0$ otherwise.  To state our results, and to facilitate our discussion of the existing literature, let us make the following definition.

\begin{defn}[Hypothesis $H_n(\delta)$]  Let $n\geq 2$ be an integer and $\delta\geq 0$.  We call {\rm Hypothesis $H_n(\delta)$} the statement that for any number field $K$, for any cuspidal automorphic representation $\pi$ of $\GL_n(\A)$ with unitary central character, and for any place $v$ of $K$, one has
\begin{equation}\label{hyp H bounds}
m(\pi,v)\leq\delta.
\end{equation}
\end{defn}

Jacquet and Shalika \cite{JS} showed\footnote{Many of the results in the literature are stated for unramified places only.  This is the case for \cite{JS}, \cite{KS}, \cite{KSh}, \cite{LRS2}, and \cite{Maki} mentioned here.  Often, however, one can prove the same numerical bounds for the ramified places with slightly more work.  For instance Rudnick-Sarnak \cite[Appendix]{RS} extend the Jacquet-Shalika bounds to ramified places, and Müller-Speh \cite{MS} do the same for the bounds of Luo-Rudnick-Sarnak.} $H_n(1/2)$ for any $n\geq 2$.   Often, however, hypothesis $H_n(1/2)$ falls just short of what is needed for concrete applications, a situation reminiscent of the subconvexity problem in the theory of $L$-functions.  It was thus a major breakthrough when Luo-Rudnick-Sarnak \cite{LRS2} showed that for any $n\geq 2$ hypothesis $H_n(\delta)$ holds for some $\delta=\delta_n<1/2$.    Their method gives the numerical value of $\delta=1/2-1/(n^2+1)$.  Owing to the existence of proven cases of functoriality in low rank, there are certain small values of $n$ for which $H_n(\delta)$ holds for a smaller $\delta$; for example Kim-Shahidi \cite{KSh} prove $H_2(1/9)$.

Over the years, various methods in analytic number theory have been developed which, for a fixed $\delta$, establish the bounds \eqref{hyp H bounds} for $K=\Q$, and possibly for $K$ an imaginary quadratic field, but not for others.  This is due to the presence of an infinite unit group for such fields.  For example, the work of Kim-Sarnak \cite{KS} establishes the bounds \eqref{hyp H bounds} for $n=2$ and any place $v$ with $\delta=7/64$, but only for the field $K=\Q$.  It has therefore been an outstanding problem to find a method robust enough to extend these results to an arbitrary number field.  A recent article of Nakasuji \cite{Maki} extends the result of Kim-Sarnak to imaginary quadratic fields.

The aim of this paper is to prove the following result, which represents an improvement over existing bounds for all fields other than $\Q$ and imaginary quadratic fields.

\begin{thm}\label{main theorem} Hypothesis $H_n(\delta_n)$ holds with $\delta_2 = 7/64$, $\delta_3 = 5/14$, $\delta_4 = 9/22$. \end{thm}

As an immediate application we obtain the following  numerical improvement for subconvexity bounds of twisted $L$-functions over number fields \cite{BH}: Let $K$ be a totally real number field, let $\pi$ be a cuspidal automorphic representation of $\GL_2(\Bbb{A})$ with unitary 
central character, and let $\chi$ be a Hecke character of conductor $\mathfrak{q}$.  Then the twisted $L$-function satisfies
\begin{displaymath}
  L(1/2, \pi \otimes \chi) \ll_{\pi, \chi_{\infty}, K, \varepsilon} \mathcal{N}(\mathfrak{q})^{\frac{1}{2} - \frac{25}{256}+\varepsilon}.
\end{displaymath}

In \cite{KS} it is shown that Theorem \ref{main theorem} follows as a consequence of the next result.

\begin{thm}\label{aux thm}
Let $\pi$ be a cuspidal automorphic representation of $\GL_n(\A)$ with unitary central character.  Assume that $L(s,\pi,{\rm sym}^2)$ converges absolutely on $\Re s>1$.  Set $m=n(n+1)/2$.  Then \eqref{hyp H bounds} holds with $\delta=\frac{1}{2}-\frac{1}{m+1}$.
\end{thm}

The rest of this introduction will serve to explain what goes into our proof of Theorem \ref{aux thm}.

\subsection{The method of Duke-Iwaniec}\label{existing}  Let $\pi$ be a cusp form on $\GL_n$ over $\Q$.  Let $\lambda_{\pi,\, {\rm sym}^2} (n)$ denote the Dirichlet coefficients of $L(s,\pi,{\rm sym}^2)$.  Fix a prime $p$.  It is an elementary exercise that the inequality \eqref{hyp H bounds} at the place $v=p$ with $\delta=(1/2)-1/(m+1)$ is equivalent to the estimate
\begin{equation}\label{show}
\lambda_{\pi,\, {\rm sym}^2} (p^\ell)\ll_\varepsilon p^{\ell(1-\frac{2}{m+1}+\varepsilon)}
\end{equation}
for arbitrarily large $\ell$. For some parameter $Q \geq 1$ let\footnote{The presence of squares in the summation condition is ubiquitous in this paper; the nice analytic properties of $L(s,\pi\otimes\chi,{\rm sym}^2)=L(s,\pi,{\rm sym}^2\times\chi^2)$ are easier to obtain than for general twists $L(s, \pi, {\rm sym}^2\times \chi)$, which have only very recently been investigated in \cite{Ta}.}   
\begin{equation}\label{Fpl}
F(p^\ell)= \sum_{\substack{Q \leq q < 2Q\\ q \not= p \text{ prime}}}\; \;  \sum_{n^2p^{-2\ell} \equiv 1 \, (q)} \lambda_{\pi,\, {\rm sym}^2}(n)g(np^{-\ell}),
\end{equation}
where $g$ is a non-negative smooth function of support in $[\frac{1}{2},2]$ satisfying $g(1)=1$.  Inverting the summation one finds
\begin{equation*}
F(p^\ell)=\sum_{p^{\ell}/2 \leq n \leq 2p^{\ell}} \lambda_{\pi,\, {\rm sym}^2}(n)f(np^{-\ell}),
\end{equation*}
where for $\gamma\in\Q^\times$ we have
\begin{equation}\label{DI test}
f(\gamma)=g(\gamma) |\{ Q\leq q< 2Q : q \text{ prime},\; q \not= p,\; \gamma^2\equiv 1\, (\text{mod } q)\}|.
\end{equation}
As $0$ has considerably more divisors than any other number, we find
\begin{equation}\label{DItrick}
|F(p^\ell)|\gg |\lambda_{\pi,\, {\rm sym}^2}(p^\ell)| \frac{Q}{\log Q} + O_\varepsilon(p^{\ell(1+\varepsilon)}).
\end{equation}
On the other hand, an upper bound for the inner sum in \eqref{Fpl} can be obtained through detecting the congruence condition by characters, inserting the functional equation for $L(s,\pi, {\rm sym}^2\times \chi^2)$ and applying Deligne's bounds on Hyper-Kloosterman sums.   In this way one can show
\begin{equation}\label{DItrick2}
F(p^\ell)\ll p^\ell+Q^{\frac{m+1}{2}+\varepsilon}.
\end{equation}
In the above error terms the absolute convergence of $L(s, \pi, \text{sym}^2)$ on $\Re s > 1$ was implicitly used.  Taking $Q= p^{2\ell/(m+1)}$ in \eqref{DItrick} and \eqref{DItrick2} gives \eqref{show}.

This argument can be adapted to the real place by setting $g$ to be the inverse Mellin transform of $L_\infty(s,\pi_\infty,{\rm sym}^2)$ and considering
\begin{equation}\label{Iwaniec} 
F(Y)=\sum_{Q\leq q < 2Q}\; \sum_{n^2 \equiv  1 \, (q)} \lambda_{\pi,\, {\rm sym}^2}(n)g(nY), \quad Y \rightarrow 0.
\end{equation}
The idea of choosing $g$ in this way to gain access to the size of the archimedean parameters is due to Iwaniec \cite{Iw}.  Alternatively, one can argue by non-vanishing of $L(s,\pi,{\rm sym}^2\times\chi^2)$ as in \cite{KS}.  It should be emphasized that the argument by non-vanishing, introduced by Luo-Rudnick-Sarnak in \cite{LRS1}, was the first treatment to successfully bound the archimedean parameters for general $\GL_n$ cusp forms, and consequently was the first to beat Selberg's $3/16$ bound on Laplacian eigenvalue of weight zero Maa{\ss} forms.

Lastly we remark that the innovation of Duke-Iwaniec \cite{DI} in the above argument was in the construction of the test function $f$ in \eqref{DI test}.  It allowed them to amplify the contribution of a single coefficient, namely $\lambda_{\pi,\, {\rm sym}^2} (p^\ell)$.  Prior to their work, one needed an $L$-series with positivite coefficients in order to drop all but one.  This in turn required the use of the Rankin-Selberg $L$-function, whose coefficients are positive but whose larger degree yields weaker bounds.  On the other hand, the Rankin-Selberg $L$-function is known to converge absolutely on $\Re s>1$ for cusp forms on $\GL_n$, making it the only tool of this sort available when dealing with an arbitrary cusp for on $\GL_n$, $n\geq 5$.

\subsection{The method of this paper}\label{ours} Let us first describe what difficulties one encounters in the situation of a general number field.   Take $\ell\geq 1$ divisible by the class number of $K$.  Let $\p$ be a prime ideal of the ring of integers $\mathcal{O}_K$ of $K$ and let $\pi$ be a generator of $\p^\ell$.  For an integral ideal $\m$ (coprime to $\p$) denote by $\mathcal{O}_K^\times\, ({\rm mod}\,\m)$ the image of the unit group $\mathcal{O}_K^\times$ in $(\mathcal{O}_K/\m)^\times$.  As an analogue of the inner sum in \eqref{Fpl} consider 
\begin{equation}\label{display13} 
 \underset{\alpha^2\pi^{-2}\, ({\rm mod}\, \m)\, \in\, \mathcal{O}_K^\times\, ({\rm mod}\,\m)}{\underset{\mathfrak{a}=(\alpha)\subseteq\mathcal{O}_K}{\sum}} \lambda_{\pi,\, {\rm sym}^2}(\mathfrak{a})g(\mathcal{N}(\mathfrak{a}\p^{-\ell})),
\end{equation}
where $g$ is as before, and $\mathcal{N}(\mathfrak{a})$ is the norm of $\mathfrak{a}$.

We immediately observe that the strength of the condition
\begin{equation}\label{bad cond}
 \alpha^2\, ({\rm mod}\, \m) \in \mathcal{O}_K^\times\, ({\rm mod}\,\m)
\end{equation}
on principal ideals $\mathfrak{a}=(\alpha)\subseteq\mathcal{O}_K$ varies with $\m$.  Indeed, for fields with an infinite unit group the image $\mathcal{O}_K^\times\, ({\rm mod}\,\m)$ can frequently be all of $(\mathcal{O}_K/\m)^\times$ (cf.\ \cite{N}), in which case the condition \eqref{bad cond} is literally empty and \eqref{DItrick} and \eqref{DItrick2} break down.  Rohrlich \cite{Ro} has a fundamental result showing that for every $\varepsilon>0$ there exists  an infinite number of square-free moduli $\m$ such that $|\mathcal{O}_K^\times\, ({\rm mod}\,\m)|\ll\mathcal{N}(\m)^\varepsilon$.  This construction was critical to the work of Luo-Rudnick-Sarnak \cite{LRS2}.  Unfortunately, the sparseness of the special moduli $\m$ given by Rohrlich's results does not allow an additional average over $\m$ as in \eqref{Fpl}. 
 
To make the above argument go through, the idea is to construct a test function on ideals that takes into account not only residue classes mod $\m$ but also certain archimedean information that compensates for the varying size of $\mathcal{O}_K^\times\, ({\rm mod}\,\m)$.  We give first a description in elementary terms for easy comparison with the previous section.  We will only need the case when $\m = \q$ is a prime ideal.  Fix a fundamental domain $\mathcal{F}$ for the action of $\mathcal{O}_K^\times$ on $K_\infty^\times=\prod_{v\mid\infty}K_v^\times$. 
Fix additionally a set of representatives $\{c\}\subseteq\mathcal{O}_K^\times$ for $\mathcal{O}_K^\times\, ({\rm mod}\,\q)$.  
Let $\mathcal{F}_c$ be the union of the translates $u \mathcal{F}$ for $u$ running over units congruent to $c$ modulo $\q$.  The condition
\begin{equation}\label{cond 1}
\alpha^2 \equiv c\,  ({\rm mod}\, \q) \quad \text{and} \quad \alpha^2 \in \mathcal{F}_c \quad \text{for some } c,
\end{equation}
is satisfied by a fraction constantly equal to
\begin{equation*}
2 \times (|\mathcal{O}_K^\times\, ({\rm mod}\,\q)|/\phi(\mathfrak{q})) \times (1/|\mathcal{O}_K^\times\, ({\rm mod}\,\q)|)= 2/\phi(\mathfrak{q})
\end{equation*}
of all principal ideals $\mathfrak{a}=(\alpha)\subseteq\mathcal{O}_K$.  More simply, one can express  \eqref{cond 1} as follows: for each principal ideal $\mathfrak{a}=(\alpha)$ let $\alpha_0$ be the unique generator such that $\alpha_0^2 \in \mathcal{F}$.  Then \eqref{cond 1} is equivalent to $\alpha_0^2 \equiv 1\pmod{\mathfrak{q}}$.

Replacing \eqref{bad cond} by \eqref{cond 1}, evaluated at $\alpha^2\pi^{-2}$, we can sum over  {\it all} primes  $\q \not = \p$  (and not just those with $|\mathcal{O}_K^\times\, ({\rm mod}\,\q)|$ small) in a dyadic interval $Q\leq\mathcal{N}(\q)<2Q$.  With this modification in place, we are able to prove the analog of estimate \eqref{DItrick} as a consequence of a suitable diophantine condition and the analog of \eqref{DItrick2} from Deligne's bounds on Hyper-Kloosterman sums.  These appear  in Section \ref{section proof} as Proposition \ref{lefthandside} and Proposition \ref{righthandside}, respectively.

We remark that the condition \eqref{cond 1} is naturally obtained by our method via an averaging operator over $S$-units.  Assume for simplicity that $\pi$ is unramified at all finite places.  Then we shall take $S=\{\q\}\cup\infty$.  Let $g_\q$ be the characteristic function of $U_\q^{(1)}=1+\q$ in $K_\q^\times$ and $g_\infty$ the characteristic function of $\mathcal{F}$ in $K_\infty^\times$.  Put $g_S=g_\q\times g_\infty$.  Then
\begin{equation*}
\sum_{u\in \mathcal{O}_S^\times}g_S(ux^2)=\begin{cases} 1,& x^2\in\mathcal{O}_S^\times(U_\q^{(1)}\times\mathcal{F});\\
										       0,&\text{else}
								 \end{cases}
\end{equation*}
is well-defined as a function on $\mathcal{O}_S^{\times}\backslash K^{\times}$.  Identifying $\mathcal{O}_S^{\times}\backslash K^{\times}$ with the group of principal ideals coprime to $\q$, and taking $g_\infty$ smooth, we thereby obtain an analytic way to capture condition \eqref{cond 1}. This convenient formalism of averaging over $S$-units will be used frequently in the present paper and allows for a clear separation between local and global properties.

Besides the above arithmetic conditions, we also need to encode some restriction of the norm into our test function, as in \eqref{display13}. When proving \eqref{hyp H bounds} for finite places, we choose $g_{\infty}$ to be some smooth approximation to the characteristic function of a ball about $1 \in K_{\infty}^{\times}$.  When proving \eqref{hyp H bounds} for an archimedean place $v \mid \infty$, we follow Iwaniec \eqref{Iwaniec} and choose $g_v$ to be the inverse Mellin transform of $L(s,\pi_v,{\rm sym}^2)$, the rest being unchanged.\\

We have compared the methods in \S\ref{existing} and the present paragraph in terms of test functions and the respective conditions on ideals they impose.  An alternative point of view is to compare them via the class of Hecke characters appearing in their spectra.  Recall that a Grö\ss encharacter $\chi$ mod $\m$ determines a unique pair $(\chi_f,\chi_\infty)$, where $\chi_f$ is a primitive character of $(\mathcal{O}_K/\m)^{\times}$, $\chi_\infty$ is a continuous homomorphism of $K_\infty^\times$ into $\C^\times$, and 
\begin{equation*}
  \chi((a))=\chi_f((a))\chi_\infty(a) \text{ for }  a\in\mathcal{O}_K,\;  (a,\m)=1.
\end{equation*}  
The character $\chi_\infty$ can be written as $\prod_{v\mid\infty}\chi_v$.  When $v=\R$ one may write $\chi_v(x)={\rm sgn}(x)^{m_v} |x|^{it_v}$, with $m_v\in\{0,1\}$ and $t_v\in\R$, and when $v=\C$ $\chi_v(z)=(z/|z|)^{m_v}|z|^{2it_v}$, with $m_v\in\Z$ and $t_v\in\R$.  Now let $C(\chi_\infty)=\prod_{v\mid\infty}(1+|m_v+it_v|)^{\deg(v)}$, where $\deg(v)=[K_v:\R]$.  The analytic conductor of $\chi$ is then $C(\chi)=C(\chi_\infty)\mathcal{N}(\m)$. This is the proper measure of complexity of a Grö\ss en\-character.  In our setting, the goal is to define a test function whose Mellin transform is supported on a set of about $X$ Grö\ss encharacters, each of conductor $X$.

This is manifestly not the case for the condition \eqref{bad cond}.  Indeed if one expands \eqref{bad cond} one obtains characters $\chi^2$ where $\chi$ is of the form $\chi (\mathfrak{a})=\omega(\mathfrak{a})\mathcal{N}(\mathfrak{a})^s$.  Here $\omega$ denotes a   ray class character of conductor dividing $\m$, so that $\omega_{\infty} = 1$ and $C(\omega) = \mathcal{N}(\m)$.  The number of such $\omega$ relative to $\mathcal{N}(\m)$ depends on the size of $\mathcal{O}_K^\times\, ({\rm mod}\,\m)$.  By contrast, the Grö\ss en\-characters obtained by expanding a smooth version of condition \eqref{cond 1} are not necessarily of finite order.  Written $\chi^2$ with $\chi(\mathfrak{a})=\omega(\mathfrak{a})\mathcal{N}(\mathfrak{a})^s$, the characters $\omega$ that contribute to the Fourier expansion in an essential way satisfy $C(\omega_{\infty}) \ll \mathcal{N}(\q)^{\varepsilon}$, whereas  $\omega_f\in \widehat{(\mathcal{O}_K/\q)^{\times}}$ may be taken arbitrary.  Thus we obtain at least $\mathcal{N}(\q)$ characters of analytic conductor at most $\mathcal{N}(\q)^{1+\varepsilon}$. This aspect of our work will be explored more in depth in \S\ref{coda}.

\subsection{Acknowledgements}  We would like to thank the Universit\'e de Nancy,  the Institute for Advanced Study and the Centre Interfacultaire Bernoulli for their hospitality and excellent working conditions.  We would like to thank P.\ Sarnak and P.\ Michel for several illuminating discussions on the subject.  We thank G.\ Henniart,  H.\ Kim and F.\ Shahidi for help in understanding some of the subtleties at ramified places, and Z.\ Rudnick for pointing out the preprint \cite{Ta}. Finally, the problem of extending some of the known techniques towards bounding Fourier coefficients to general number fields was suggested to the second named author by A.\ Venkatesh, who also emphasized the importance of the use of infinite order Hecke characters.
 
\section{Some $\GL_1$ preliminaries}\label{prelim}

Let $K$ be a number field of degree $d $ over $\Q$.  We write $r$ for the number of inequivalent archimedean embeddings of $K$.  Let $\mathcal{O}_K$ be the ring of integers of $K$ and $\mathfrak{d}$ the different.  Let $P_K=\mathcal{O}_K^\times\backslash K^\times$ be the group of all principal fractional ideals of $K$.  Let $h$ be the class number of $K$.  For an integral ideal $\mathfrak{a}$ we denote its norm by $\mathcal{N}(\mathfrak{a})$ and put $\phi(\mathfrak{a})=|(\mathcal{O}_K/\mathfrak{a})^{\times}|$.

For each place $v$ of $K$ let $|\cdot |_v$ be the normalized $v$-adic absolute value.  In particular, if $v=\Bbb{C}$ then $|\cdot|_v$ is the square of the modulus.  Write $K_v$ for the completion of $K$ with respect to $|\cdot |_v$. If $v=\mathfrak{p}$ let $\mathcal{O}_\mathfrak{p}$ be the ring of integers of $K_\p$, $\mathfrak{d}_\p$ the local different, and $\varpi_\p$ a uniformizer.  Let $U_\p=  U_\p^{(0)} = \mathcal{O}_{\mathfrak{p}}^\times$ and   $U_\p^{(r)}=1+\p^r$ for $r \geq 1$. We  set $U_v=\{\pm 1\}$ if $K_v = \R$ and $U_v= U(1)$ if $K_v=\C$.  In all cases, $U_v$ is the maximal compact subgroup of $K_v^{\times}$.

For any locally compact abelian group $A$ let $\widehat{A}$ be the group of characters.  These are the continuous homomorphisms into $U(1)$.  It will sometimes be convenient to write a character $\chi\in\widehat{K_v^\times}$ as $\chi(x)=|x|_v^{it}\eta(x)$, where $t\in\R$ and $\eta\in\widehat{U_v}$.  When $v=\R$ we write $\eta_m(x)={\rm sgn}(x)^m$ for $m\in\{0,1\}$.  When $v=\C$ we write $\eta_m(z)=(z/|z|)^m$ for $m\in\Z$.  When $v=\p$ the real number $t$ is determined only up to an integer multiple of $2\pi/\log\mathcal{N}(\p)$.  We say $\eta\in\widehat{U_\p}$ has degree $m\geq 1$, and write $\deg(\eta)=m$, if $\eta$ is trivial on $U_\p^{(m)}$ but not on $U_\p^{(m-1)}$.  In all cases, $\chi(x)=|x|_v^{it}\eta(x)$ is said to be unramified when $\eta=1$. 

Let $S$ be a finite set of places containing all infinite places.  Set $K_S^\times=\prod_{v\in S}K_v^\times$ endowed with its norm $|\cdot |_S=\prod_{v\in S}|\cdot |_v$.  Denote by $\mathcal{O}_S$ the ring of $S$-integers.  Let $P_K(S)$ be the group of principal fractional ideals prime to $S$.  Then $P_K(S)$ can be identified with $\mathcal{O}_S^\times\backslash K^\times$ by sending the ideal $(\gamma) \in P_K(S)$ to the orbit $\mathcal{O}_S^{\times}.\gamma$.  The inverse map sends the $\mathcal{O}_S^\times$-orbit $\mathfrak{o}$ to the ideal $(\gamma)$, where $\gamma$ is any element in $\mathfrak{o}$ with $v_{\mathfrak{p}}(\gamma) = 0$ for all finite $\mathfrak{p} \in S$.  Under this identification the norm $\mathcal{N}(\gamma)$ of the ideal $(\gamma)$ is $|\gamma|_S$.  We denote by $\Delta: \R_+ \hookrightarrow K_\infty^\times$ the map $t\mapsto\prod_{v\mid\infty} t^{1/d}$.

Let $\A$ be the adele ring of $K$.  Let $\I= \A^{\times}$ be the group of ideles of $K$.  Put $|\cdot |_\A$ for the idelic norm and let $\I^1$ be the closed subgroup of $\I$ consisting of ideles of norm 1. Put $\mathscr{C}=K^\times\backslash\I$ and $\mathscr{C}^1=K^\times\backslash\I^1$, the latter of which is compact by Dirichlet's theorem.  We identify $\widehat{\I^1}$ (resp.\ $\widehat{\mathscr{C}^1}$) with the closed subgroup of $\widehat{\I}$ (resp.\ $\widehat{\mathscr{C}}$) consisting of those characters trivial on $\Delta(\R_+)\subset\I$.  We have $\widehat{\I}\simeq\widehat{\I^1}\times\R$ and
\begin{equation}\label{isomo}
\widehat{\mathscr{C}}\simeq\widehat{\mathscr{C}^1}\times\R,
\end{equation}
the correspondence $\chi\leftrightarrow (\omega,t)$ being given by $\chi(x)=\omega(x)|x|_\A^{it}$.  If $\m$ is an integral ideal of $\mathcal{O}_K$ we denote by $\widehat{\I}(\m)$ (resp.  $\widehat{\mathscr{C}}(\m)$, $\widehat{\mathscr{C}^1}(\m)$) the group of characters $\chi \in \widehat{\I}$ (resp., $\widehat{\mathscr{C}}$, $\omega \in \widehat{\mathscr{C}^1}$) of conductor dividing $\m$.

We fix a non-trivial character $\psi_v$ of $K_v$ by taking $\psi_v(x)=\exp(2\pi i x)$ when $v=\R$, $\psi_v(x)=\exp(2\pi i(x+\overline{x}))$ when $v=\C$, and $\psi_\p$ an additive character trivial on $\mathfrak{d}_{\mathfrak{p}}^{-1}$ but not on $\varpi^{-1}_\p\mathfrak{d}_{\mathfrak{p}}^{-1}$ when $v=\p$ is non-archimedean.  Let $dx_v$ be the self-dual Haar measure on $K_v$.  Explicitly, $dx_v$ is Lebesque measure if $v$ is real, twice the Lebesque measure if $v$ is complex, and the unique Haar measure such that $\mathcal{O}_\mathfrak{p}$ has volume $\mathcal{N}(\mathfrak{d}_\p)^{-1/2}$ if $v=\p$ is finite.  On $K_v^{\times}$ we choose the normalized Haar measure $d^\times x_v = \zeta_v(1)dx_v/|x_v|_v$, where $\zeta_v$ is the Tate local zeta function at $v$.  We let $d^\times x$ be the measure on $\I$ that on the standard basis of open sets of $\I$ coincides with $ \prod_v d^\times x_v$.  We continue to denote by $d^\times x$ the quotient measure on $\mathscr{C}$.

Let $A$ be one of the groups $K_v^\times$, $\I$, or $\mathscr{C}$, taken with its norm $|\cdot |$ and the choice of Haar measure indicated above, which we write here as $d^\times a$.  For $\sigma \in \Bbb{R}$ let $L^1(A, \sigma) = \{ g : | \cdot |^{\sigma} g \in L^1(A)\}$.  When $g\in L^1(A, \sigma)$ we write $\widehat{g}(\sigma, \chi) = \int_A g(a) \chi(a) |a|^{\sigma}d^{\times}a$.  If $g$ is continuous and $\widehat{g}(\sigma, .) \in L^1(\widehat{A})$ the Mellin inversion formula reads
\begin{equation}\label{inversion} 
g(a) = \int_{\widehat{A}} \widehat{g}(\sigma, \chi) \chi^{-1}(a) |a|^{-\sigma} d\chi
\end{equation}
for a unique choice of Haar measure $d\chi$ on $\widehat{A}$.  When $A=K_v^\times$ and $\chi(x)=|x|_v^{it}\eta(x)$ we sometimes write $\widehat{g}(s,\eta)$ in place of $\widehat{g}(\sigma,\chi)$, where $s=\sigma+it$.  Similarly, when $A=\I$ or $\mathscr{C}$, and $\chi(x)=|x|_\A^{it}\omega(x)$, we sometimes write $\widehat{g}(s,\omega)$ in place of $\widehat{g}(\sigma,\chi)$.  The measures $d\chi$ on $\widehat{K_v^\times}$ are explicitly given by
\begin{displaymath}
c_v \sum_m \int_{(\sigma)} g(s, \eta_m) \frac{ds}{2\pi i}, \qquad \sum_{\eta \in \widehat{U_{\mathfrak{p}}}} \int_{\sigma -\frac{i\pi}{\log \mathcal{N}(\mathfrak{p})}}^{\sigma+\frac{i \pi}{\log \mathcal{N}(\mathfrak{p})}} g(s, \eta) \log\mathcal{N}(\mathfrak{p})\frac{ds}{2\pi i},
\end{displaymath}
for $v\mid\infty$ and $\mathfrak{p}$, respectively. Here $c_\R=1/2$ and $c_\C=1/(2\pi)$.  The corresponding measure on $\widehat{\mathscr{C}}\simeq\widehat{\mathscr{C}^1}\times\R$ is the product of the counting measure on the first factor and $c_K/2\pi$ times Lebesgue measure on the second factor,  where $c_K^{-1} =\underset{s=1}{\rm Res}\; \zeta_K(s)$ (see \cite[VII.6. Prop.\ 12]{We}). 

For fixed $c\in\R$ let $\mathfrak{z}(c)$ be the vector space of continuous complex-valued functions $g$ on $\I$ such that the $K^\times$-invariant function $G(x)=\sum_{\gamma\in K^\times} g(\gamma x)$ converges absolutely and uniformly on compacta in $\I$, and in addition  $g \in L^1(\I, \sigma)$, $\widehat{G}(\sigma,\cdot ) \in L^1(\widehat{\mathscr{C}})$ for all $\sigma > c$.  The unfolding technique shows
\begin{equation*} 
\widehat{g}(\sigma,\chi)=\int_{\Bbb{I}}g(x)\chi(x)|x|_\A^{\sigma}d^{\times}x=\int_{\mathscr{C}}G(x)\chi(x)|x|_\A^{\sigma}d^{\times}x=\widehat{G}(\sigma,\chi).
\end{equation*}
Clearly $G\in L^1(\mathscr{C},\sigma)$ for all $\sigma>c$.  We deduce that for $g\in\mathfrak{z}(c)$ the function $G$ satisfies the criteria under which \eqref{inversion} holds.  Thus
\begin{equation}\label{unfolding}
G(x)=\int_{\widehat{\mathscr{C}}}   \widehat{G}(\sigma,\chi) \chi^{-1}(x) |x|_\A^{-\sigma} d\chi=\int_{\widehat{\mathscr{C}}}\widehat{g}(\sigma,\chi) \chi^{-1}(x) |x|_\A^{-\sigma} d\chi 
\end{equation}
for $\sigma>c$.
\section{Symmetric square $L$-functions}\label{sym2 L}

Let $n\geq 2$.  Let $\pi$ be a cuspidal automorphic representation of $\GL_n(\A)$, with unitary central character $\omega_\pi$.  Fix an identification $\pi\simeq \otimes_v'\pi_v$, where $\pi_v$ is an irreducible unitary representation of ${\rm GL}_n(K_v)$.  Denote by $\widetilde{\pi}=\otimes_v \widetilde{\pi}_v$ the contragredient representation of $\pi$.  We shall always normalize $\pi$ so that $\omega_\pi$ is trivial on $\Delta(\R_+)\subset\I$.  We may take an arbitrary $\pi$ into this form by twisting it by $|\det |^{it}$ for an appropriate $t\in\R$.

Let $k$ be a local or global field.  Let $\GL_n(\C)\times W_k$ be the $L$-group of $G=\GL_n$, where $\GL_n$ is viewed as an algebraic group over $k$.  Here $W_k$ is the Weil-Deligne group of $k$.  For $K$ a global field and $K_v$ the completion of $K$ at the place $v$, there is a natural map $\theta_v: W_{K_v}\rightarrow W_K$.  Thus if $\rho$ is a finite dimensional complex representation of $\GL_n(\C)\times W_K$ then there is an associated collection $\{\rho_v\}$ of finite dimensional complex representations of $\GL_n(\C)\times W_{K_v}$, each given by composition with $Id\times\theta_v$.

We return to the case where $K$ is a number field.  Let $\pi\simeq\otimes_v\pi_v$ be as above.  Let $\rho$ be a finite-dimensional representation of $\GL_n(\C)\times W_K$.  To this data Langlands has attached an Euler product $\Lambda (s,\pi,\rho)=\prod_v L(s,\pi_v,\rho_v)$.  We shall be interested in character twists of the symmetric square representation ${\rm sym}^2\times\chi^2:\GL_n(\C)\times W_K\rightarrow \GL_m(\C)$, where $m = n(n+1)/2$.

The local symmetric square $L$-function $L(s,\pi_v,{\rm sym}^2)$ at $v$ is defined by Shahidi in \cite{Sh81}, for $v$ infinite, and \cite{Sh0}, for $\p$ finite.  For $\p$ finite, it is of the form $P(\mathcal{N}(\p)^{-s})^{-1}$ for a polynomial with complex coefficients and constant term 1.  For all $v$ if $\pi_v$ is tempered then $L(s,\pi_v,{\rm sym}^2)$ is holomorphic on $\Re s>0$ (see \cite{Sh85} for $v$ archimedean).  For $\pi_v$ non-tempered, then writing it as a Langlands quotient as in the introduction we have the factorisation (see, for example, \cite[page 30]{Sh92})
\begin{equation*}
L(s,\pi_v,{\rm sym}^2)=\prod_{j=1}^rL(s+2\sigma_j,\tau_j,{\rm sym}^2)\prod_{i<j}L(s+\sigma_i+\sigma_j,\tau_i\times\tau_j).
\end{equation*}
Since local Rankin-Selberg $L$-functions are holomorphic on $\Re s>0$ for tempered pairs (see, for example, \cite[Appendix]{RS}), we deduce that in all cases $L(s,\pi_v,{\rm sym}^2)$ is holomorphic on $\Re s>2m(\pi,v)$, the real number $m(\pi,v)$ being defined in the introduction.

For $\p$ finite we expand $L(s,\pi_\p,{\rm sym}^2)$ into a Dirichlet series, obtaining
\begin{equation*}
L(s,\pi_\p,{\rm sym}^2)=\sum_{r\geq 1}\lambda_{\pi,\, {\rm sym}^2} (\p^r)\mathcal{N}(\p)^{-rs},
\end{equation*}
for some coefficients $\lambda_{\pi,\, {\rm sym}^2} (\p^r)$, satisfying $\lambda_{\pi, \, {\rm sym}^2} (1)=1$.  This series converges absolutely on $\Re s>2m(\pi,\p)$.  Since $\widetilde{\pi}_\p\simeq\overline{\pi}_\p$, the coefficients of $L(s,\widetilde{\pi}_\p,{\rm sym}^2)$ are simply $\overline{\lambda}_{\pi,\, {\rm sym}^2} (\p^r)$.  

For $\Re s > m(\pi,v)$ put $Z(s,\chi_v,\pi_v,{\rm sym}^2)=L(2s, \pi_v,{\rm sym}^2\times\chi_v^2)$ for $\chi_v$ unramified and zero otherwise.  As a function of $\chi_v$ it is in $L^1(\widehat{K^{\times}_v})$ for all $\Re s > m(\pi,v)$.  For $\sigma$ in this range we  define the inverse Mellin transform 
\begin{equation}\label{coeff fun}
\lambda_v(x)=\int_{\widehat{K_v^\times}}Z(\sigma,\chi_v,\pi_v,{\rm sym}^2)\chi_v^{-1}(x)|x|_v^{-\sigma}d\chi_v.
\end{equation}
Thus $\lambda_v$ is continuous and $U_v$-invariant. More explicitly, for $v \mid \infty$ one has
\begin{equation*}
\lambda_v(x)= c_v\int_{(\sigma)} L(2s,\pi_v,{\rm sym}^2)|x|_v^{-s}\frac{ds}{2\pi i},
\end{equation*}
where $c_\R=1/2$ and $c_\C=1/(2\pi)$, while for $v=\p$ one has
\begin{align}\label{explicit coeff}
\lambda_\p(x)&=\int_{\sigma-i\pi/\log\mathcal{N}(\p)}^{\sigma+i\pi/\log\mathcal{N}(\p)}L(2s,\pi_\p,{\rm sym}^2)|x|_\p^{-s}\log\mathcal{N}(\p)\frac{ds}{2\pi i}\nonumber\\
&=\begin{cases}
\lambda_{\pi,\, {\rm sym}^2}(\p^r),& \text{for } v_\p(x)=2r, r\geq 0;\\
0,&\text{otherwise}.
\end{cases}
\end{align}
Observe that $\lambda_\p(\varpi_\p^{2r}) = \lambda_{\pi, {\rm sym}^2}(\p^r)$.  In all cases, shifting the contour  we find that
\begin{equation}\label{blow-up}
\lambda_v(x) \begin{cases} \ll_\sigma |x|_v^{-\sigma},\qquad\text{as}\; x\rightarrow 0,\quad \forall\; \sigma > m(\pi,v);\\
\ll_{A} |x|_v^{-A}, \qquad \text{as}\; |x| \rightarrow \infty.\end{cases}
\end{equation}
Bounding the blow-up rate of $\lambda_v$ at zero is therefore equivalent to bounding $m(\pi,v)$.  This observation seems to have been first used by Iwaniec in \cite{Iw}.  It allows us to treat all places in a uniform way.  Finally, it is easy to see that $\widehat{\lambda_v}(\sigma, \chi) = Z(\sigma, \chi, \pi_v, {\rm sym}^2)$ for  $\sigma > m(\pi,v)$. 

Let $S$ be any finite set of places of $K$ containing all infinite places.   Denote by $L^S(s,\pi,{\rm sym}^2)$ the product of $L(s,\pi_\p,{\rm sym}^2)$ over all $\p\notin S$.  We have
\begin{equation*}
L^S(s,\pi, {\rm sym}^2)=\underset{(\mathfrak{a},S)=1}{\underset{\mathfrak{a}\subseteq\mathcal{O}_K}{\sum}}\lambda_{\pi,\, {\rm sym}^2} (\mathfrak{a})\mathcal{N}(\mathfrak{a})^{-s},
\end{equation*}
where $\lambda_{\pi,\, {\rm sym}^2}(\mathfrak{a})=\prod_{\p^r||\mathfrak{a}}\lambda_{\pi,\, {\rm sym}^2} (\p^r)$ for $(\mathfrak{a},S)=1$.  This series converges absolutely for $\Re s > 3/2$ by the bounds $H_n(1/2)$ of Jacquet-Shalika.

Let $B_{\pi}$ be the set of places at which $\pi$ is ramified, together with all infinite places. Assume   that $S$ contains $B_\pi$.  Then it was  proven in \cite{BG} (see also \cite{Ta}) and then again in \cite{K1} by the Langlands-Shahidi method that  the function $L^S(s,\pi, {\rm sym}^2)$ admits a meromorphic continuation to all of $\C$.  In fact, if $\pi\not\simeq\widetilde{\pi}$ then $L^S(s,\pi, {\rm sym}^2)$ is entire; whereas if $\pi\simeq\widetilde{\pi}$, the only possible poles are finite in number\footnote{Of course it is believed that there are no poles in the critical strip but it is enough for us to know that the number of such poles is $O_{\pi}(1)$.} and located within the critical strip.  This follows immediately from  \cite[Theorem 4.1]{Ta}.  Moreover, the work \cite{GS} assures that away from possible poles, the function $L^S(s,\pi,{\rm sym}^2)$ is of moderate growth on vertical lines.

Shahidi \cite{Sh1} has shown that the functional equation
\begin{equation}\label{fnl eq 1}
L^S(s,\pi,{\rm sym}^2)=\gamma_S(s,\pi,{\rm sym}^2)L^S(1-s,\widetilde{\pi},{\rm sym}^2)
\end{equation}
holds for all $s\in \C$, where $\gamma_S(s,\pi,{\rm sym}^2)=\prod_{v\in S}\gamma (s,\pi_v,{\rm sym}^2)$ satisfies $\gamma(s,\pi_v,{\rm sym}^2)=\epsilon(s,\pi_v,{\rm sym}^2)L(1-s,\widetilde{\pi}_v,{\rm sym}^2)/L(s,\pi_v,{\rm sym}^2)$.  Moreover
\begin{enumerate}

\item\label{unram fin} for finite $\p\notin B_{\pi}$ one has $\epsilon(s,\pi_\p,{\rm sym}^2)=1$; thus $\gamma(s,\pi_\p,{\rm sym}^2)=P_\p(\mathcal{N}(\p)^{-s})/Q_\p(\mathcal{N}(\p)^{-(1-s)})$ for   polynomials $P_\p, Q_\p$ of degree $m$ such that $P_\p(0)=Q_\p(0)=1$.  Moreover $Q_\p(\mathcal{N}(\p)^{-(1-s)})\neq 0$ for $\Re s\leq 0$;

\item\label{ram fin} for finite $\p\in B_{\pi}$ the function $\gamma(s,\pi_\p,{\rm sym}^2)$ is a rational function in $\mathcal{N}(\p)^{-s}$ that is pole free for $\Re s\leq 0$;

\item\label{inf gam} for infinite $v\mid\infty $ the function $\gamma(s,\pi_v,{\rm sym}^2)$ is a meromorphic function of moderate growth in vertical strips (away from its poles), and is holomorphic on $\Re s<0$, pole free on $\Re s=0$.\footnote{The pole free regions given in \eqref{ram fin} and \eqref{inf gam} above can be improved by the work of Luo-Rudnick-Sarnak.  This improvement will not be needed as an input to our method.} 

\end{enumerate}

For later purposes we record several reformulations of the above facts upon replacing $\pi$ by $\pi\otimes\chi$.  We have $L^S(s,\pi\otimes\chi,{\rm sym}^2)=L^S(s,\pi,{\rm sym}^2\times\chi^2)$, where by $\chi^2$ on the right hand side we intend the unique character of $W_K$ associated to $\chi^2$ via the (dual of the) homeomorphism $\mathscr{C}\overset{\sim}{\longrightarrow}W_K^{ab}$.  One has the factorization $L^S(s,\pi,{\rm sym}^2\times\chi^2)=\prod_{\p\notin S}L(s,\pi_\p,{\rm sym}^2\times\chi_\p^2)$.  The functional equation \eqref{fnl eq 1} becomes
\begin{equation}\label{fnl eq 2}
L^S(s,\pi,{\rm sym}^2\times\chi^2)=\gamma_S(s,\pi,{\rm sym}^2\times\chi^2)L^S(1-s,\widetilde{\pi},{\rm sym}^2\times\chi^{-2}).
\end{equation}
If $\pi\otimes\chi$ is not self-dual then $L^S(s,\pi,{\rm sym}^2\times\chi^2)$ is holomorphic.  On the other hand, if $\pi\otimes\chi\simeq\widetilde{\pi}\otimes\chi^{-1}$ then by equating central characters we deduce that $\chi^n=\overline{\omega}_\pi$.

Fix a prime $\q\notin B_\pi$.  It will be useful to quantify how many $\chi\in\widehat{\mathscr{C}}(\mathfrak{q})$ can satisfy $\chi^n=\overline{\omega}_\pi$.  Note that the conductor $\m$ of $\omega_\pi$ has support in $B_\pi$.  Since $\q\notin B_\pi$, we see that if $\m\neq {\bf 1}$ then no such $\chi$ can verify $\chi^n=\overline{\omega}_\pi$.

\begin{lemma}\label{quantify} Let $\q$ be a prime ideal of $K$.  Let $\xi\in\widehat{\mathscr{C}}$ be fixed, of conductor ${\bf 1}$.  Then the number of $\chi\in\widehat{\mathscr{C}}(\q)$ such that $\chi^n=\xi$ is $O(1)$, where the implied constant depends only on $K$ and $n$.
\end{lemma}

{\bf Proof.}  For the proof we use the language of Grö\ss en\-characters.  Assume the conductor of $\chi$ is $\q$, the case where the conductor is equal to ${\bf 1}$ being similar.  Recall that $\chi$ mod $\q$ determines (uniquely, up to multiplication by a class group character) a pair $(\chi_f,\chi_\infty)$, where $\chi_f$ is a primitive character of $(\mathcal{O}_K/\q)^\times$, $\chi_\infty=\prod_{v|\infty}\chi_v$ is a character of $K_\infty^\times$, and $\chi ((a))=\chi_f(a)\chi_\infty(a)$ for all $a\in\mathcal{O}_K$ prime to $\q$.  Similarly $\xi$ determines a character $\xi_\infty$ of $K_\infty^\times$.  We may therefore equate $\chi_v=\xi_v$ for every $v|\infty$ and $\chi_f^n=1$ and then count the number of solutions in each equation individually.  If $v=\R$ there are at most two $\chi_v$ such that $\chi_v^n=\xi_v$, depending on the parity of $\xi_v$ and $n$. If $v=\C$ there is at most one $\chi_v$ such that $\chi_v^n=\xi_v$.  Since the group $(\mathcal{O}_K/\q)^{\times}$ is cyclic of order $\phi(\q)$, there are at most $(n,\phi(\q))\leq n$ choices of $\chi_f$.\qed\\

Let $B_{\pi,K}$ be the set of finite places at  which $K$ is ramified, together with all places in $B_{\pi}$.  Let $\q\notin B_{\pi,K}$ be a prime and take $\chi\in\widehat{\mathscr{C}}(\mathfrak{q})$.  Then
\begin{equation}\label{gamma twist 1}
\gamma(s,\pi_\q,{\rm sym}^2\times\chi_\q^2)=L(1-s,\widetilde{\pi}_\q,{\rm sym}^2\times\chi_\q^{-2})/L(s,\pi_\q,{\rm sym}^2\times\chi_\q^2)
\end{equation}
for $\chi$ of conductor ${\bf 1}$, and
\begin{equation}\label{gamma twist 2}
\gamma(s,\pi_\q,{\rm sym}^2\times\chi_\q^2)=\epsilon(s,\pi_\q,{\rm sym}^2\times\chi_\q^2)=\mathcal{N}(\q)^{-ms}\tau(\chi_\q^2)^m
\end{equation}
for $\chi$ of conductor $\q$.  Here $\tau(\chi_\q)=\sum_\varepsilon \chi_\q(\varepsilon)\psi_\q(\varpi_\q^{-1}\varepsilon)$ is the Gau{\ss} sum, where $\varepsilon$ runs through a set of representatives of $U_\q/U_\q^{(1)}$.

\section{Local and $S$-adic computations}\label{local section}

We make the following general conventions that remain valid for the rest of the paper. As in the previous section, we let $n\geq 2$ and fix $\pi$ a cuspidal automorphic representation of $\GL_n(\A)$, with unitary central character $\omega_\pi$.

For $v \in B_{\pi}$ let $g_v$ be a smooth, $U_v$-invariant function on $K_v^\times$ that is either of compact support or equal to $\lambda_v$.   Next fix a prime ideal $\q\notin B_{\pi,K}$ and let $g_\q$ be the characteristic function on $U_\q^{(1)}$.  Write $S = \{\q\}\cup B_\pi$.  Let $T\subset S$ be the subset of places where the function $g_v$ is a genuine ``test function", that is, $T=\{\q\} \cup \{v \in B_{\pi} : g_v \neq \lambda_v\}$.  For each $v\in S$ define
\begin{equation}\label{define}
g_v^*(x)=\int_{\widehat{K_v^\times}} \widehat{g_v}(1/2-\sigma,\chi^{-1}) \gamma(1-2\sigma, \pi_v, \text{sym}^2\times\chi^2) \chi^{-1}(x)|x|_v^{-\sigma} d\chi
\end{equation}
where $\sigma>1/2$.  By Mellin inversion
\begin{equation}\label{fstar dual}
\widehat{g_v^*}(\sigma,\chi)=\widehat{g_v}(1/2-\sigma,\chi^{-1}) \gamma(1-2\sigma, \pi_v, \text{sym}^2\times\chi^2).  
\end{equation}

\smallskip

\begin{lemma}\label{dual supp}  {\rm (1)} We have
\begin{equation*}
\begin{cases} 
		      g^*_\q(x) \ll_\varepsilon \phi(\q)^{-1}(|x|_\q^{-\frac{1}{2}-\varepsilon} + |x|_\q^{-\frac{1}{4}+\frac{1}{4m}}),& |x|_\q\leq \mathcal{N}(\q)^{2m};\\
                      g^*_\q(x) =0,& |x|_\q>\mathcal{N}(\q)^{2m}.
\end{cases} 
\end{equation*}
{\rm (2)}  For $v \in B_{\pi}$,  $A \geq 1$, and $0<\varepsilon<1/2$ we have $g^*_v(x)\ll_{v, \varepsilon} \min( |x|_v^{-\frac{1}{2}-\varepsilon}, |x|_v^{-A})$.
 \end{lemma}

{\bf Proof.} (1) One computes easily that $\widehat{g_\q}=\phi(\q)^{-1}{\bf 1}_{\deg(\chi)\leq 1}$.  Using the explicit formula for $\gamma(1-2s,\pi_\q,{\rm sym}^2\times\chi^2)$ in \eqref{gamma twist 1} and \eqref{gamma twist 2} we find $g_\q^*(x)=\phi(\q)^{-1}(A(x)+B(x))$, where
\begin{equation*}
A(x)=\int_{\sigma-i\pi/\log\mathcal{N}(\q)}^{\sigma+i\pi/\log\mathcal{N}(\q)}\frac{L(2s,\widetilde{\pi}_\q,{\rm sym}^2)}{L(1-2s,\pi_\q,{\rm sym}^2)}|x|_\q^{-s}\log\mathcal{N}(\q)\frac{ds}{2\pi i}
\end{equation*}
and
\begin{equation*}
B(x)=\sum_{\deg (\eta)= 1}\tau(\eta^2)^m\overline{\eta}(x/|x|_\q)\int_{\sigma-i\pi/\log\mathcal{N}(\q)}^{\sigma+i\pi/\log\mathcal{N}(\q)}\mathcal{N}(\q)^{m(2s-1)}|x|_\q^{-s}\log\mathcal{N}(\q)\frac{ds}{2\pi i}.
\end{equation*}

The integrand $L(2s,\widetilde{\pi}_\q,{\rm sym}^2)/L(1-2s,\pi_\q,{\rm sym}^2)$ in $A(x)$ is described in \eqref{unram fin} of \S\ref{sym2 L}.  A direct calculation then shows that $A(x)=0$ if $|x|_\q>\mathcal{N}(\q)^{2m}$.  Otherwise, we shift the line of integration to $\sigma = 1/2+\varepsilon$ for any $\varepsilon>0$ (which is admissible by the Jacquet-Shalika bounds $H_n(1/2)$) for $0 < |x|_\q \leq 1$ and to some very large number for $1<|x|_\q\leq\mathcal{N}(\q)^{2m}$, and estimate trivially.

The integral in $B(x)$ is non-zero if and only if $|x|_\q=\mathcal{N}(\q)^{2m}$ in which case it is $\mathcal{N}(\q)^{-m}$.  Thus we have
\begin{equation*}
B(x)=|x|_\q^{-1/2}\sum_{\deg (\eta)= 1}\tau(\eta^2)^m\overline{\eta}(x/|x|_\q),\qquad (|x|_\q=\mathcal{N}(\q)^{2m}).
\end{equation*}
Now for $y\in U_\q$ we have
\begin{equation*}
\frac{1}{\phi(\q)}\bigg((-1)^m+\sum_{\deg (\eta)= 1}\tau(\eta^2)^m\overline{\eta}(y)\bigg)=\sum_{(y_1\cdots y_m)^2=y}\psi_\q(\varpi_\q^{-1}(y_1+\cdots +y_m)).
\end{equation*}
The latter convolution sum vanishes if $y$ is not equal to $y_0^2$ for some $y_0\in U_\q$, and otherwise if $y=y_0^2$ it is   $Kl(y_0)+Kl(-y_0)$, where
\begin{equation*}
Kl(y)=\sum_{y_1\cdots y_m=y}\psi_\q(\varpi_\q^{-1}(y_1+\cdots +y_m))
\end{equation*}
is the degree-$m$ Hyper-Kloosterman sum.  Deligne \cite[p.\ 219]{De1} has shown $Kl(y)\ll \mathcal{N}(\q)^{\frac{m-1}{2}}$, so that $B(x)\ll |x|_\q^{-\frac{1}{2}}\mathcal{N}(\q)^{\frac{m+1}{2}}=|x|_\q^{-\frac{1}{4}+\frac{1}{4m}}$ for $|x|_\q=\mathcal{N}(\q)^{2m}$.

(2) By properties \eqref{ram fin} and \eqref{inf gam}  in  \S\ref{sym2 L} we can shift the contour in \eqref{define} to $\Re s = 1/2+\varepsilon$ resp.\ to $\Re s = A$ getting the desired bounds. \qed\\

The following useful result of Bruggeman-Miatello \cite[Lemma 8.1]{BM} states, essentially, that in descending local estimates to global ones, one only loses on a logarithmic scale.

\begin{lemma}\label{bruggemia} Let $a, b \in \Bbb{R}$, $a+b > 0$. Let $g : K_{\infty}^\times \rightarrow \Bbb{C}$ be a function satisfying $|g(x)| \leq \prod_{v \mid \infty} \min(|x_v|_v^a, |x_v|_v^{-b})$. Then
\begin{displaymath}
  \sum_{u \in \mathcal{O}_K^\times} |g(ux)| \ll_{ a, b} \left(1+|\log |x|_{\infty}|^{r-1}\right) \min (|x|_{\infty}^a, |x|_{\infty}^{-b}).
\end{displaymath} 
\end{lemma}

With the notation and assumptions as in the beginning of this section put $g_S=\prod_{v\in S}g_v$ and $g^*_S=\prod_{v\in S}g^*_v$.  Next we write
\begin{equation}\label{defGS}
G_S(x)=\sum_{u\in\mathcal{O}_S^\times}g_S(ux)\quad\text{and}\quad G^*_S(x)=\sum_{u\in\mathcal{O}_S^\times}g_S^*(ux).
\end{equation}

\begin{prop}\label{S-adic}  We have
\begin{equation*}
\begin{cases}
G^*_S(x)\ll_{\varepsilon, \pi} \mathcal{N}(\q)^{\frac{m-1}{2}+\varepsilon} |x|_S^{-\frac{1}{2}}, &|x|_S\geq 1;\\
 G^*_S(x)\ll_{\varepsilon, A, \pi} |x|_S^{-A}, & |x|_S\geq \mathcal{N}(\q)^{2m+\varepsilon}
 \end{cases}
\end{equation*}
for all $A\geq 1$ and $\varepsilon > 0$. 
\end{prop}

{\bf Proof.} Let $S':= S \setminus \{\mathfrak{q}\}=B_\pi$.  By  Lemma \ref{dual supp} we have, for $|x|_S\geq 1$,
\begin{equation*}
g^*_S(x)\ll_{\varepsilon, \pi} \phi(\q)^{-1}|x|_S^{-\frac{1}{2}-\varepsilon}\mathcal{N}(\mathfrak{q})^{2m(\frac{1}{4} + \frac{1}{4m} +\varepsilon})\prod_{v \in S'} \min (1,  |x_v|_v^{-A}),
\end{equation*}
and $g_S^*(x)=0$ if $|x|_{\mathfrak{q}} > \mathcal{N}(\mathfrak{q})^{2m}$. 
We fix a set of representatives of $\mathcal{O}_K^{\times}\backslash \mathcal{O}_S^{\times}$. By Lemma \ref{bruggemia} we obtain $G_S^{\ast}(x)\ll_{ \varepsilon, A, \pi} \mathcal{N}(\q)^{\frac{m-1}{2}+\varepsilon}|x|_S^{-\frac{1}{2}} K$ where 
\begin{displaymath}
  K 
 =  \sum_{\substack{u \in \mathcal{O}_S^{\times}\backslash \mathcal{O}_K^{\times}\\v_{\mathfrak{q}}(ux_{\mathfrak{q}}) \geq -2m} } \prod_{\mathfrak{p} \in S'}\min(1, |ux_\mathfrak{p}|_{\mathfrak{p}}^{-A})    \min\left(|ux_{\infty}|_{\infty}^{-\varepsilon} , |ux_{\infty}|_{\infty}^{-A}\right) 
\end{displaymath}
for any $A\geq 0$, $\varepsilon > 0$. 
We can majorize $K$ as follows:
\begin{displaymath}
\begin{split}
  K \leq & \sum_{\ell_{\mathfrak{q}} \geq -2m - v_{\mathfrak{q}}(x_\mathfrak{q})}\sum_{\substack{\ell_{\mathfrak{p}} \in \Bbb{Z}\\ \mathfrak{p} \in S'}} \prod_{\mathfrak{p} \in S'}\min\left(1, (\mathcal{N}\mathfrak{p}^{-\ell_{\mathfrak{p}}} |x_\mathfrak{p}|_{\mathfrak{p}})^{-A}\right)\\
  & \times   
  \min\Bigl(\Bigl(| x_{\infty}|_{\infty}\mathcal{N}(\mathfrak{q})^{\ell_{\mathfrak{q}}} \prod_{\mathfrak{p} \in S'} \mathcal{N}(\mathfrak{p})^{\ell_\mathfrak{p}}\Bigr)^{-\varepsilon}, \Bigl(| x_{\infty}|_{\infty}\mathcal{N}(\mathfrak{q})^{\ell_{\mathfrak{q}}} \prod_{\mathfrak{p} \in S'} \mathcal{N}(\mathfrak{p})^{\ell_\mathfrak{p}}\Bigr)^{-A}\Bigr).
  \end{split}
\end{displaymath}
Let $r$ be the number of finite primes in $S'$.  Inductively it is now easy to see that 
\begin{displaymath}
\begin{split}
  K& \ll_{\varepsilon, A}  \sum_{\ell_{\mathfrak{q}} \geq -2m - v_{\mathfrak{q}}(x_\mathfrak{q})} \min\Bigl(\Bigl(| x_{S'}|_{S'}\mathcal{N}(\mathfrak{q})^{\ell_{\mathfrak{q}}}  \Bigr)^{-(r+1)\varepsilon}, \Bigl(| x_{S'}|_{S'}\mathcal{N}(\mathfrak{q})^{\ell_{\mathfrak{q}}}  \Bigr)^{-A+r\varepsilon}\Bigr)\\
 & =  \sum_{\ell_{\mathfrak{q}} \geq -2m} \min\Bigl(\Bigl(| x_{S}|_{S}\mathcal{N}(\mathfrak{q})^{\ell_{\mathfrak{q}}}  \Bigr)^{-(r+1)\varepsilon}, \Bigl(| x_{S}|_{S}\mathcal{N}(\mathfrak{q})^{\ell_{\mathfrak{q}}}  \Bigr)^{-A+r\varepsilon}\Bigr). 
\end{split}
\end{displaymath}
The bounds in the proposition are now obvious. \qed
 
\section{Voronoi formula}

The goal of this section is to prove the summation formula in Proposition \ref{Vor}.  We begin by a technical lemma.

Let $\lambda^S = \prod_{v\not\in S} \lambda_v$.  Write $g=\lambda^S \times g_S$ and $g^*=\overline{\lambda^S} \times g^*_S$ as complex valued functions on $\Bbb{I}$.

\begin{lemma}\label{function lemma} The functions $g$ and $g^*$ lie in the space $\mathfrak{z}(c)$ for any $c \geq 1$.\end{lemma}

{\bf Proof.} It is clear that $g$ and $g^*$ are continuous. To prove that $g$ and $g^*$ are in $L^1(\I,\sigma)$, for any $\sigma >1$, we first observe that for $\mathfrak{p}\notin S$,
\begin{align*}
\|  \lambda_{\mathfrak{p}}\|_{L^1(K_{\mathfrak{p}}^{\times}, \sigma)} &= \int_{K_\p^\times}\lambda_\p(x)|x|_\p^\sigma d^\times x=\sum_{r\geq 0}\lambda_{\pi,{\rm sym}^2}(\p^r)\mathcal{N}(\p)^{-2r\sigma}\\
&=1+O_\varepsilon(\sum_{r\geq 1}\mathcal{N}(\p)^{-2r(\sigma-m(\pi,\p)-\varepsilon)})=1+O_\varepsilon (\mathcal{N}(\p)^{2(m(\pi,\p)-\sigma+\varepsilon)})
\end{align*}
and
\begin{displaymath}
\begin{split}
  \| \widehat{\lambda}_\p(\sigma, \cdot) \|_{L^1(\widehat{K_\p^{\times}} )} & = \int_{\sigma-\pi/\log\mathcal{N}(\p)}^{\sigma+\pi/\log\mathcal{N}(\p)} | L(2s,\pi_\p,{\rm sym}^2)| \log \mathcal{N}(\p) \frac{|ds|}{2\pi}\\
  & = 1 + O_\varepsilon(\mathcal{N}(\p)^{2(m(\pi,\p)-\sigma+\varepsilon)}). 
 \end{split} 
\end{displaymath}
In particular, if $\sigma-m(\pi,v)>1/2+\varepsilon$ for every $v$, then the infinite products
\begin{equation*}
\|  g \|_{L^1(\Bbb{I}, \sigma)} = \prod_v \|  g_{v}\|_{L^1(K_{v}^{\times}, \sigma)},\qquad \| \widehat{g}(\sigma, \cdot) \|_{L^1(\widehat{\Bbb{I}} )} = \prod_v \| \widehat{g}_{v}(\sigma, \cdot)\|_{L^1(\widehat{K_{v}^{\times}})}
\end{equation*}
converge.  The same is true when $g$ is replaced by $g^*$.  By the Jacquet-Shalika bounds $H_n(1/2)$ we may take any $\sigma>1$.

Next we prove the absolute and  uniform convergence of $G$ on compacta $C$.  By the $K^{\times}$-invariance of $G$ we may assume that $v_{\p}(x) \leq 0$ for all $\p$ and all $x \in C$.  Recall that $\supp(g_\p)=\mathcal{O}_\p$ and $g_\p(x)\ll |x|_\p^{-\frac{1}{2}-\varepsilon}$ for $\p\notin S$, while $g_v(x) \ll \min(|x|_v^{-\frac{1}{2}-\varepsilon}, |x|_v^{-2})$ for $v\in S$.  Letting $U=\prod_{\p\not\in S}\mathcal{O}_\p\times K_S^\times$ we find
\begin{displaymath}
\sum_{\gamma \in K^{\times}} |g(\gamma x)|\ll_C \sum_{\gamma\in K^{\times} \cap x^{-1}U}\; \prod_{v \in S} \min(1, |\gamma x_v|_v^{-\frac{3}{2}+\varepsilon} ).
\end{displaymath}
The sum over $\gamma\in K^{\times} \cap x^{-1}U$ can be regrouped along $\mathcal{O}_S^\times$-cosets $(\alpha)\in P_K(S)=\mathcal{O}_S^\times\backslash K^\times$ such that $v_\p(\alpha)\geq-v_\p(x) \geq 0$ for all $\p\notin S$. The resulting $S$-unit sum is
\begin{equation*}
  \sum_{u\in \mathcal{O}_S^{\times}}\prod_{v \in S} \min(1, | u \alpha x_v|_v^{-\frac{3}{2}+\varepsilon}).
\end{equation*}
The same argument as in Proposition \ref{S-adic} can be used to bound the preceding expression by $|\alpha x_S|_S^{\varepsilon-3/2}$.  We can complete the $\alpha$-sum to a sum over all integral ideals obtaining  $\sum |g(\gamma x)|\ll_C  \zeta_K(3/2-\varepsilon)$ for $x \in C$.  The same argument works for $G^{\ast}$. 

The fact that $\int_{\widehat{\mathscr{C}}} |\widehat{G}(\chi)| d\chi < \infty$ and $\int_{\widehat{\mathscr{C}}} |\widehat{G}(\chi)| d\chi < \infty$ follows easily from the decay properties of $g_S$ and $g_S^{\ast}$ after the same manipulations as in \eqref{expansion}, \eqref{exp1} below.
\qed\\

We have
\begin{equation}\label{Mellin f}
\widehat{g}(\sigma,\chi)=\begin{cases}
    L^S(2\sigma, \pi, \text{sym}^2\times \chi^2) \widehat{g_S}(\sigma,\chi), & \chi\in\widehat{\I}(\q); \\
    0, & \text{else}
    \end{cases}  
\end{equation}
and
\begin{equation}\label{gstarhat}
\widehat{g^*}(\sigma,\chi)=\begin{cases}
    L^S(2\sigma, \widetilde{\pi}, \text{sym}^2\times \chi^{-2}) \widehat{g_S^*}(\sigma,\chi), & \chi\in\widehat{\I}(\q); \\
    0, & \text{else}
    \end{cases}  
\end{equation}
for $\sigma > 1$.  The existence of the left hand sides of \eqref{Mellin f} and \eqref{gstarhat}  is guaranteed by Lemma \ref{function lemma}. Let
\begin{equation}\label{def of G and G*}
G(x)=\sum_{\gamma\in K^\times}g(\gamma x)\quad\text{and}\quad G^*(x)=\sum_{\gamma\in K^\times}g^*(\gamma x),
\end{equation}
as functions on $\mathscr{C}$.  We deduce from Lemma \ref{function lemma} and \eqref{unfolding} that
 \begin{equation}\label{F}
G(x) = \int_{\widehat{\mathscr{C}}} \widehat{g}(\sigma, \chi)\chi^{-1}(x)|x|_\A^{-\sigma} d\chi,\,  G^*(x)=\int_{\widehat{\mathscr{C}}}\widehat{g^*}(\sigma, \chi)  \chi^{-1}(x)|x|_\A^{-\sigma} d\chi
 \end{equation}
for any $\sigma>1$.

\begin{prop}[Voronoi summation]\label{Vor} Let the notation and assumptions be as described in Section \ref{local section}.  Then
\begin{equation*}
G(x^2) =  |x|_\A^{-1}R + |x|_\A^{-1} G^*(1/x^2)
\end{equation*}
for $x \in \mathscr{C}$, where
\begin{equation}\label{def of R}
R=c_K\underset{\omega^n=\overline{\omega}_\pi}{\underset{\omega \in \widehat{\mathscr{C}^1}(\q)}{\sum}}  \overline{\omega}^2(x)\sum_{\rho} \underset{s=\rho}{\rm Res}\   |x|_\A^{1-2s}\widehat{g_T}(s,\omega)  L^T (2s, \pi, {\rm sym}^2\times \omega^2),
\end{equation}
with  $c_K^{-1} =   \underset{s=1}{\rm Res}\; \zeta_K(s)$ and the sum over $\rho$ running over all poles of  $L^T(2s,\pi,{\rm sym}^2\times\omega^2)$. The (possibly empty) sum over $\omega$ and $\rho$ is finite. 
\end{prop} 
 
\textbf{Proof.} We apply the decomposition \eqref{isomo} to \eqref{F} to obtain
\begin{equation}\label{expansion}
G(x^2)=c_K\sum_{\omega\in\widehat{\mathscr{C}^1}}\overline{\omega}^2(x)\int_{(3/2)} \widehat{g}(s,\omega)|x|_\A^{-2s} \frac{ds}{2\pi i}.
\end{equation}
By \eqref{Mellin f}, we have  
\begin{equation}\label{exp1}
\widehat{g}(s,\omega) =
    L^S(2s, \pi, \text{sym}^2\times \omega^2) \widehat{g_S}(s,\omega) =  L^T(2s, \pi, \text{sym}^2\times \omega^2) \widehat{g_T}(s,\omega),
\end{equation}
if $\omega\in\widehat{\mathscr{C}^1}(\q)$ and $\widehat{g}(s, \omega) = 0$ otherwise. 
Hence we can restrict the sum over  $\omega$ to the set $\widehat{\mathscr{C}^1}(\q)$. For each such $\omega$ we shift the contour to $\Re s = -1$.  This is admissible by the rapid decay of the infinite components of $\widehat{g}(s, \omega)$ along vertical lines.  We apply the functional equation \eqref{fnl eq 2} and change variables $s\mapsto 1/2-s$, $\omega\mapsto\overline{\omega}$.  In this way we obtain
\begin{align*}
& G(x^2) = |x|_\A^{-1}R+c_K |x|_\A^{-1} \sum_{\omega \in \widehat{\mathscr{C}^1}(\q)} \overline{\omega}^2(1/x)\\
 &\times  \int_{(3/2)} L^S(2s,\widetilde{\pi}, \text{sym}^2\times \overline{\omega}^2)\widehat{g_S}(1/2-s,\overline{\omega}) \gamma_S(1-2s, \pi, \text{sym}^2\times\omega^2)|1/x|_\A^{-2s} \frac{ds}{2\pi i}.
 \end{align*}
By the known properties of the poles of $L(s,\pi,{\rm sym}^2)$ and Lemma \ref{quantify} the sum defining $R$ is finite.  Applying \eqref{isomo} in the other direction together with \eqref{fstar dual} we find
\begin{align*}
G(x^2)&=|x|_\A^{-1}R+|x|_\A^{-1}\int_{\widehat{\mathscr{C}}} L^S(2\sigma,\widetilde{\pi}, \text{sym}^2\times \chi^{-2})\widehat{g^*_S}(\sigma,\chi)\chi(1/x^2)^{-1}|1/x^2|_\A^{-\sigma}d\chi.
\end{align*}
From \eqref{F} and \eqref{gstarhat} we deduce that the last term above equals $|x|_\A^{-1}G^*(1/x^2)$.

\section{Proof of Theorem \ref{aux thm}}\label{section proof}

We prove Theorem \ref{aux thm} by combining the Voronoi summation formula given in Proposition \ref{Vor} with the $S$-adic estimates in Proposition \ref{S-adic}.  Recall the hypothesis of Theorem \ref{aux thm} that $L(s,\pi,{\rm sym}^2)$ converges absolutely on $\Re s>1$.

We first establish Theorem \ref{aux thm} for $v_0 =\mathfrak{p}_0 \not\in B_{\pi}$ finite and then indicate the necessary changes for an infinite or ramified place at the end of the section.   We keep the notation and assumptions as in Section \ref{local section}.

Fix $\q \not\in B_{\pi, K}$ distinct from $\p_0$.  For each $v\in B_{\pi}$ we take $g_v$ of compact support.  Thus $T=S=\{\q\}\cup B_\pi $.  Moreover, for every $v\in B_{\pi}$ we require that $g_v$ be non-negative and satisfy $g_v(1)=1$.  Let $x\in\Bbb{I}$ be such that $x_v=1$ for $v \not= \p_0$ and $x_{\p_0}=\varpi_{\p_0}^\ell$ for a positive integer $\ell$ divisible by class number $h$ of $K$; clearly $|x|_\A = \mathcal{N}(\p_0)^{-\ell}$.

We first write the functions $G(x^2)$ and $G^*(1/x^2)$ as smooth sums of Dirichlet series coefficients.  To this end, recall that $\lambda_\p$ is $U_\p$-invariant for each $\p \not \in S$ and $P_K(S)=\mathcal{O}_S^\times\backslash K^\times$; hence
\begin{equation}\label{changelater}
G(x)=\sum_{(\gamma)\in P_K(S)}\lambda^S(\gamma x^S)G_S(\gamma x_S)
\end{equation}
with $G_S$ as in \eqref{defGS}. Inputting our choice of $x$ this gives
\begin{displaymath}
G(x^2)=\sum_{(\gamma) \in P_K(S)}  \bigg(\prod_{\p\notin S\cup \{\p_0\}}\lambda_\p(\gamma)\bigg)\lambda_{\p_0}(\gamma\varpi_{\p_0}^{2\ell})G_S(\gamma). 
\end{displaymath}
Now $\lambda_\p$ for $\p\not\in S$ is supported on $\mathcal{O}_\p$.  The summation over $(\gamma)$ therefore runs through the semi-group of fractional principal ideals prime to $S$ of the form $\mathfrak{a} \p_0^{-\ell}$ where $\mathfrak{a}$ is an integral principal ideal prime to $S$.  Changing variables and inserting \eqref{explicit coeff} we obtain
\begin{equation}\label{simp sum 1}
G(x^2)=\sum_{\substack{\mathfrak{a} \subseteq\mathcal{O}_K\\ \mathfrak{a} \in P_K(S)}}  \lambda_{\pi,\, {\rm sym}^2}(\mathfrak{a})G_S(\mathfrak{a}^2 \p_0^{-2\ell}).
\end{equation}
Likewise
\begin{equation}\label{simp sum 2}
G^*(1/x^2)=\sum_{\substack{\mathfrak{a} \subseteq\mathcal{O}_K\\ \mathfrak{a} \in P_K(S)}}  \overline{\lambda}_{\pi,\, {\rm sym}^2}(\mathfrak{a})G_S^*(\mathfrak{a}^2 \p_0^{2\ell}).
\end{equation}
Recall that if $\mathfrak{c}=(\gamma)\in P_K(S)$ then $G_S(\mathfrak{c})$ simply means $G_S(\gamma)$.  

 As we now want to vary $\q$, our notation will henceforth explicitly reflect the dependence on all quantities on $\q$.  Thus we write $S_\q=\{\q\}\cup B_\pi$, and correspondingly we define $g_{S_\q}$ and $G_{S_\q}$.   Moreover we now write $g (x;\q)=\lambda^{S_\q}\times g_{S_\q}$, $G(x;\q)=\sum_{\gamma\in K^\times} g(\gamma x;\q)$, and similarly for $g^*(x;\q)$ and $G^*(x;\q)$.

Fix a parameter $Q\gg_{K,\pi} \log \mathcal{N}(\p_0)$, to be chosen later in terms of $\mathcal{N}(\p_0)$.  Let $\mathcal{Q}=\{\q \text{ prime}: Q\leq\mathcal{N}(\q)<2Q,\; \q\neq \p_0,\; \q\notin B_{\pi,K}\}$ and put $F(x)=\sum_{\q\in\mathcal{Q}}G(x;\q)$.

\begin{prop}\label{lefthandside} There is a constant $c > 0$ depending only on $K$ such that
\begin{displaymath}
|F(x^2)|\geq c |\lambda_{\pi,\, {\rm sym}^2}(\p_0^\ell)| \frac{Q}{\log Q} + O_\varepsilon \left(\mathcal{N}(\p_0)^{(1+\varepsilon)\ell}\right).
\end{displaymath}
\end{prop}

\textbf{Proof.} Put $S':=S_\q\setminus \{\q\}=B_\pi$.  We switch the order of summation in $F(x^2)$ and use \eqref{simp sum 1} to obtain
\begin{equation*}
\begin{split}
F(x^2)&=\sum_{\substack{\mathfrak{a}\subseteq\mathcal{O}_K\\ \mathfrak{a}\in P_K(S')}}\lambda_{\pi,\, {\rm sym}^2}(\mathfrak{a})\sum_{\q\in\mathcal{Q}: (\q,\mathfrak{a})=1}G_{S_\q}(\mathfrak{a}^2\p_0^{-2\ell})\\
&  =\lambda_{\pi,\, {\rm sym}^2}(\p_0^\ell)\sum_{\q\in\mathcal{Q}: (\q,\mathfrak{a})=1}G_{S_\q}(\textbf{1})+E,
\end{split}
\end{equation*}
where $E$ is the sum over all $\mathfrak{a}\neq\p_0^\ell$. As $G_{S_\q}(\textbf{1})\geq 1$,  we have
\begin{equation*}
\sum_{\q\in\mathcal{Q}: (\q,\mathfrak{a})=1}G_{S_\q}(\textbf{1})\geq \sum_{\q\in\mathcal{Q}:(\q,\mathfrak{a})=1}1 \gg_K \frac{Q}{\log Q}.
\end{equation*}

As usual, for $x\in K_{S'}^\times$ let $G_{S'}(x)$ denote the average of $g_{S'}(ux)$ as $u$ ranges over $\mathcal{O}_{S'}^\times$.  Note that we have $G_{S'}(x)=0$ for $|x|_{S'}$ outside of a closed bounded interval in $(0,\infty)$ depending only on the support conditions on $g_v$, $v\in B_\pi$.

Now suppose that for a principal fractional ideal $\mathfrak{c} \in P_K(S')$, $\mathfrak{c}\neq {\bf 1}$, we have
\begin{equation}\label{claim 2}
\sum_{\q\in\mathcal{Q}}G_{S_\q}(\mathfrak{c})\ll_\varepsilon H(\mathfrak{c})^{\varepsilon}G_{S'}(\mathfrak{c}), 
\end{equation}
where $H(\mathfrak{c})= \max(\mathcal{N}(\mathfrak{a}), \mathcal{N}(\mathfrak{b}))$ if $\mathfrak{c} = \mathfrak{a}/\mathfrak{b}$ for coprime integral ideals $\mathfrak{a}, \mathfrak{b}$.  We deduce from \eqref{claim 2} that
\begin{equation*}
E\ll_\varepsilon \mathcal{N}(\p_0)^\varepsilon \sum_{\mathcal{N}(\mathfrak{a})\ll\mathcal{N}(\p_0)^\ell} |\lambda_{\pi,\, {\rm sym}^2}(\mathfrak{a})|\mathcal{N}(\mathfrak{a})^\varepsilon.
\end{equation*}
The absolute convergence of $L(s,\pi,{\rm sym}^2)$ on $\Re s > 1$ implies $E\ll_\varepsilon\mathcal{N}(\p_0)^{(1+\varepsilon)\ell}$, yielding the proposition.

To substantiate the claim \eqref{claim 2} let $\mathfrak{c}=(\gamma)$ for $\gamma\notin\mathcal{O}_K^\times$ and $(\gamma,S')=1$.  Then
\begin{displaymath}
  \sum_{\mathfrak{q} \in \mathcal{Q}} G_{S_\q}(\gamma) = 
 \sum_{\mathfrak{q} \in \mathcal{Q}} \sum_{ u \in \mathcal{O}_{S_\q}^{\times} } g_{S_\q}(u\gamma) =  \sum_{\mathfrak{q} \in \mathcal{Q}} \sum_{\substack{u \in \mathcal{O}_{S'}^\times\\ u\gamma \equiv 1 \, (\text{mod }\q)}} g_{S'}(u\gamma) \leq G_{S'}(\gamma) W(\gamma),
\end{displaymath}
where
\begin{displaymath}
W(\gamma) = \underset{ \substack{u\in \mathcal{O}_{S'}^{\times}\\ u\gamma \in \text{supp}(g_{S'})}}{\sup}| \{\q \in \mathcal{Q} : u\gamma \equiv 1 \, (\text{mod }\q)\}|.
\end{displaymath}
To estimate $W(\gamma)$, recall that if $\gamma=\alpha/\beta$ for $\alpha,\beta\in\mathcal{O}_K$, the notation $u\gamma\equiv 1\pmod\q$ means $\q|(u\alpha-\beta)$.  Our assumption that $\mathfrak{c}\neq (1)$ translates to the principal integral ideal $(u\alpha-\beta)$ being non-zero.  Hence
\begin{equation*}
W(\gamma)\ll_\varepsilon    \underset{ \substack{u\in \mathcal{O}_{S'}^{\times}\\ u\gamma \in \text{supp}(g_{S'})}}{\sup}   \mathcal{N}(u\alpha-\beta)^\varepsilon.
\end{equation*}
The $\sup$ runs over $S'$-units $u$ satisfying $u\alpha-\beta\in\beta.\Omega$ where $\Omega$ is the neighborhood of $0\in K_{S'}$ given by $\supp(g_{S'})-1$.  Thus $\mathcal{N}(u\alpha-\beta)\ll_{g_{S'}} \mathcal{N}(\beta)$ for all such $u$.  Now if $\mathfrak{c}=\mathfrak{a}\mathfrak{b}^{-1}$ for relatively prime integral ideals $\mathfrak{a},\mathfrak{b}$, then $\mathcal{N}(\beta)\leq \mathcal{N}(\mathfrak{b})^h\leq H(\mathfrak{c})^h$.  This shows that $W(\gamma)\ll_{\varepsilon,K, \pi} H(\mathfrak{c})^\varepsilon$, as desired.\qed
 
\begin{prop}\label{righthandside} We have $F(x^2)\ll_\varepsilon  \mathcal{N}(\p_0)^\ell+Q^{\frac{m+1}{2}+\varepsilon}$.\end{prop}

\textbf{Proof.} It suffices to prove that $G(x^2;\q)\ll_\varepsilon  \mathcal{N}(\q)^{-1}\mathcal{N}(\p_0)^\ell+\mathcal{N}(\q)^{\frac{m-1}{2}+\varepsilon}$.  From Proposition \ref{Vor} we have $G(x^2;\q) = |x|_\A^{-1} R_\q+ |x|_\A^{-1} G^*(1/x^2;\q)$.

We have $\widehat{g_{S_\q}}(\rho,\omega)\ll\phi(\q)^{-1}$ for $\Re\rho\in [0,1/2]$.  By Lemma \ref{quantify}, the number of $\omega\in\widehat{\mathscr{C}^1}(\q)$ such that $\omega^n=\overline{\omega}_\pi$ is $O(1)$.  Moreover, for any given $\omega$ such that $\pi\otimes\omega$ is self-dual, the number of poles of $L^{S_\q}(s,\pi\otimes\omega,{\rm sym}^2)$ is $O_n(1)$.  We conclude that $R_\q\ll \mathcal{N}(\q)^{-1}$.

From \eqref{simp sum 2} and Proposition \ref{S-adic} we deduce
\begin{equation*}
G^*(1/x^2;\q)\ll_\varepsilon \mathcal{N}(\q)^{\frac{m-1}{2}+\varepsilon}\mathcal{N}(\p_0)^{-\ell}\sum_{\mathcal{N}(\mathfrak{a})\ll_\varepsilon \mathcal{N}(\q)^{m+\varepsilon}\mathcal{N}(\p_0)^{-\ell}}|\lambda_{\pi,\, {\rm sym}^2}(\mathfrak{a})|\mathcal{N}(\mathfrak{a})^{-1}.
\end{equation*}
This last sum is $\ll_\varepsilon \mathcal{N}(\q)^\varepsilon$, by the absolute convergence of $L(s,\pi,{\rm sym}^2)$ on $\Re s>1$.  From this we deduce $|x|^{-1}F^*(1/x^2)\ll_\varepsilon \mathcal{N}(\q)^{\frac{m-1}{2}+\varepsilon}$, as desired.\qed\\

{\bf Proof of Theorem \ref{aux thm}.}  Propositions \ref{lefthandside} and \ref{righthandside} combine to give
\begin{equation*}
\lambda_{\pi,\, {\rm sym}^2}(\p_0^\ell) \ll Q^{\varepsilon-1} \mathcal{N}(\p_0)^{(1+\varepsilon)\ell} + Q^{\frac{m-1}{2}+\varepsilon}.
\end{equation*}
Choosing $Q=\mathcal{N}(\p_0)^{\frac{2\ell}{m+1}}$ we obtain $\lambda_{\pi,\, {\rm sym}^2}(\p_0^\ell)\ll_\varepsilon \mathcal{N}(\p_0)^{\ell(1-\frac{2}{m+1}+\varepsilon)}$, which is
\begin{equation}\label{finite bound}
\lambda_{\p_0}(x)\ll_\sigma |x|_{\p_0}^{-\sigma}, \qquad x=\varpi_{\p_0}^{2\ell},\quad \forall\; \sigma >\frac{1}{2}-\frac{1}{m+1}.
\end{equation}
Letting $\ell\rightarrow\infty$ (along multiples of $h$) we obtain from \eqref{blow-up} the desired bound $m(\pi,\p_0)\leq (1/2)-1/(m+1)$.\\

The case $v_0 = \mathfrak{p}_0 \in B_{\pi}$ is almost identical, with one notable exception: for the local function at $\p_0$ we take $g_{\mathfrak{p}_0} = \lambda_{\mathfrak{p}_0}$.  The rest of the argument requires only notational changes.  We leave the details to the reader, and instead sketch the necessary modifications to the argument at archimedean $v_0$.

When $v_0 \mid \infty$ we take $g_{v_0}=\lambda_{v_0}$.  As before $g_\q$ is the characteristic function of $U_\q^{(1)}$, and for every finite prime $\p\in B_\pi$, $g_\p$ is the characteristic function of $U_\p$.  For each $v\in \infty\setminus\{v_0\}$ we may choose $g_v$ to satisfy $g_v(1)=1$ and have small enough support so that $g_S(u)=0$ if $u\in\mathcal{O}_S^\times$ is not a root of unity.  Let $x \in \Bbb{I}$ satisfy
\begin{equation}\label{infty x}
x_v = 1\; \text{for } v \not= v_0\, \text{and } x_{v_0} = Y \; \text{for a small parameter }  0<Y<1.
\end{equation}
The parameter $Y$ will tend to 0; it plays the role of $|\varpi_0^{\ell}|_{\p_0} =  \mathcal{N}(\p)^{-\ell}$ for $\ell$ large.  As in \eqref{simp sum 1}, for this choice of $x$ we have
\begin{equation}\label{infty G}
G(x^2;\q) = \sum_{\substack{\mathfrak{a} = (\gamma) \in P_K(S)\\ \mathfrak{a} \subseteq \mathcal{O}_K}} \lambda_{\pi, {\rm sym}^2}(\mathfrak{a}) G_{S_\q}(\gamma^2 x_S^2).
\end{equation}
When $\gamma=1$ we obtain $G_S(x_S^2;\q)=w_K\lambda_{v_0}(Y^2)$, where $w_K$ is the number of roots of unity in $K$.

The argument of Proposition \ref{lefthandside} can be adapted to the current situation.  As before, one obtains the inequality \eqref{claim 2}.  In this case $G_{S'}$ no longer vanishes for $|x|_{S'}$ large enough, but applying the local bound \eqref{blow-up} at the place $v_0$ and Lemma \ref{bruggemia} one shows that $G_{S'}(x)\ll_{\varepsilon, A} \min(|x|_{S'}^{-1/2-\varepsilon}, |x|_{S'}^{-A})$.   This is enough for the same argument to go through.  We find that
\begin{equation}\label{low F}
|F(x^2)|\geq c |\lambda_{v_0}(Y^2)| \frac{Q}{\log Q} + O_\varepsilon(Y^{-1-\varepsilon}),
\end{equation}
the main term coming from $\mathfrak{a} = (1)$. Next 
\begin{equation*}
G^*(1/x^2;\q)=\sum_{\substack{\mathfrak{a} \subseteq\mathcal{O}_K\\ \mathfrak{a} \in P_K(S_\q)}}  \overline{\lambda}_{\pi,\, {\rm sym}^2}(\mathfrak{a})G_{S_\q}^*(\gamma^2/x_S^2).
\end{equation*}
We apply Proposition \ref{S-adic} unchanged, and arguing as in the proof of Proposition \ref{righthandside} one finds $F(x^2) \ll Y^{-1}+Q^{\frac{m+1}{2}+\varepsilon}$.

Choosing $Q = Y^{-2/(m+1)}$ as before we find, similarly to \eqref{finite bound}, that
\begin{equation*}
\lambda_{v_0}(x)\ll_\sigma |x|_{v_0}^{-\sigma}, \qquad x=Y^2,\; 0<Y< 1,\quad \forall\; \sigma >\frac{1}{2}-\frac{1}{m+1}.
\end{equation*}
The bounds $m(\pi,v_0)\leq \frac{1}{2}-\frac{1}{m+1}$ now follow from \eqref{blow-up}.\qed

\section{Coda}\label{coda}

Our goal in this section is expository; we thank the referee for suggesting to include it.  We give an alternative proof of Theorem \ref{aux thm} in the case of $K$ a real quadratic field.  This proof is similar in spirit to the method of Luo-Rudnick-Sarnak in \cite{LRS1}, \cite{LRS2}, where bounds towards Ramanujan are deduced from the non-vanishing of character twists of certain $L$-functions.

Let $\mu,\nu$ denote the two real embeddings of $K$.  We prove Theorem \ref{aux thm} at the place $\nu$.
 
Let $\q\notin B_{\pi, K}$ be a prime ideal and let $g_\q$ be the characteristic function of $U_\q^{(1)}$.  For finite primes $\p\in B_\pi$ let $g_\p$ be the characteristic function of $U_\p$.  Let $g_\mu\in C_c^\infty(\R^\times)$ be even, non-negative, and satisfy $g_\mu(1)=1$.  For a given $\beta\in (0,1/2)$ let $g_\nu\in C^\infty(\R^\times)$ be even, non-negative, satisfy $g_\nu(y)=0$ for $|y|\geq 2$, and $g_\nu(y)= |y|^{-\beta}$ for  $0 < |y| < 1/2$.  As $g_\mu$ and $g_\nu$ are even, their Mellin transforms are zero on characters of the form ${\rm sgn}(x)|x|^{it}$.  Note that $\widehat{g}_\nu(\sigma, |\cdot |^{it})=\widehat{g}_\nu(s)$, viewed as a meromorphic function in $s=\sigma+it$, has a pole at $s=\beta$ of residue 1, while $\widehat{g}_\mu(s)$ is entire; both are of rapid decrease in vertical strips.  

Let $S=\{\q\}\cup B_\pi$.  We write $g_S = \prod_{v\in S}g_v$.   For $v\notin S$ put $g_v=\lambda_v$.  Consider $G(x)$ as in \eqref{def of G and G*}.  Since $g_\p$ is $U_\p$-invariant for all finite $\p\neq\q$ the function $\widehat{G}(\sigma,\chi)$ is supported on $\widehat{\mathscr{C}(\q)}$.

For every fixed $\chi\in\widehat{\mathscr{C}}(\q)$, written $\chi=|\cdot |_\A^{it}\omega$ where $\omega=\prod_v\omega_v\in\widehat{\mathscr{C}^1}$ as in \eqref{isomo}, the function of a complex variable $s=\sigma+it$ given by $\widehat{g}(\sigma,\chi)=\widehat{g}(s,\omega)=\prod_v\widehat{g}_v(s,\omega_v)=L^S(2s,\pi, {\rm sym}^2\times\omega^2)\widehat{g_S}(s,\omega_S)$ is holomorphic except possibly at 
\begin{enumerate}

\item\label{first pole} a finite set, of cardinality at most $O_\pi(1)$, of simple poles in the critical strip $0\leq 2\Re(s)\leq 1$ when $\omega$ is such that $\pi\otimes\omega\simeq\widetilde{\pi}\otimes\overline{\omega}$, coming from the $L$-factor $L^S(2s,\pi, {\rm sym}^2\times\omega^2)$,

\item\label{sec pole} and at $s=\rho_\omega=\beta-ir_\nu$, where $\omega_\nu=|\cdot |^{ir_\nu}$, coming from $\widehat{g_\nu}(\rho_\omega,\omega_\nu)=\widehat{g_\nu}(\beta) = \infty$.

\end{enumerate}
(By adjusting $\beta$ by an $\varepsilon$ amount, these two sets of poles can be assumed disjoint.)  We emphasize that even when $\omega$ satisfies $\pi\otimes\omega\simeq\widetilde{\pi}\otimes\overline{\omega}$, it is not necessarily the case that $L^S(2s,\pi, {\rm sym}^2\times\omega^2)$ has a pole at $s=1/2$.  Moreover at the point described in \eqref{sec pole} it could very well happen that $L^S(2s,\pi, {\rm sym}^2\times\omega^2)$ vanishes, killing the pole.  Our argument will show that the latter cannot take place too often.

Consider the expansion \eqref{expansion} for some $x \in \Bbb{I}$ to be chosen momentarily.  For every fixed $\omega\in\widehat{\mathscr{C}^1}(\q)$ we shift the vertical contour across the critical strip, picking up the possible poles enumerated in \eqref{first pole} and \eqref{sec pole}.  Let $T=S\setminus\{\nu\} $. Following the proof of Proposition \ref{Vor} we obtain $G(x^2)= |x|_\A^{-1}R+|x|_\A^{-2\beta}L(x)+|x|_\A^{-1}G^*(1/x^2)$, where $R$ is defined in \eqref{def of R} and 
\begin{equation}\label{L}
L(x):=c_K\sum_{\substack{\chi \in \widehat{\mathscr{C}}(\mathfrak{q})\\\ \chi_\nu = 1}} L^S(2\beta , \pi, {\rm sym}^2 \times \chi^2) \widehat{g_T}(\beta,\chi_T)|x|_\A^{-it_\mu}.
\end{equation}
In the above, we have written $\chi_\mu(x)=|x|^{it_\mu}{\rm sgn}^{m_\mu}$.  To see \eqref{L} note for each $\omega\in\widehat{\mathscr{C}(\q)}$ the point $s=\rho_\omega$ of case \eqref{sec pole} contributes
\begin{equation*}
c_K|x|^{-2\rho_\omega}L^S(2\rho_\omega,\pi,{\rm sym}^2\times\omega^2)\widehat{g_T}(\rho_\omega,\omega_T).
\end{equation*}
Next observe that $\rho_\omega$ can be descibed as the unique point such that the character $\chi$ defined by $|\cdot |_\A^{\rho_\omega}\omega=|\cdot |_\A^\beta \chi$ satisfies $\chi_\nu=1$.  The triviality of $\chi_\nu$ implies $r_\nu=-t_\mu/2$, and we may then write $\rho_\omega=\beta-ir_\nu=\beta+it_\mu/2$ to arrive at \eqref{L}.

For $\varepsilon>0$ let $X_\varepsilon (\q)$ denote the set of Hecke characters $\chi\in\widehat{\mathscr{C}}(\q)$ such that $ \chi_\nu = 1$, ${\rm Cond}(\chi_\mu) \leq \mathcal{N}(\q)^\varepsilon$, and ${\rm deg}(\chi_\q) \leq 1$.  Then
\begin{equation}\label{newL}
L(x)= \frac{c_K}{\phi(\q)} \sum_{\chi\in X_\varepsilon (\q)} L^S(2\beta , \pi, {\rm sym}^2 \times \chi^2)c(\chi, x)+ O_{A,\varepsilon}(\mathcal{N}(\q)^{-A}),
\end{equation}
where $c(\chi, x)=\widehat{g_\mu}(\beta,\chi_\mu)|x|_\A^{-it_\mu}\ll 1$.  Is easy to verify that $|X_\varepsilon(\q)|\asymp\mathcal{N}(\q)^{1+\varepsilon}$.  Note furthermore that since we do not additionally require that $\chi_\mu=1$, the characters in $X_\varepsilon(\q)$ are not necessarily of finite order.

As before, we now reintroduce the dependence on $\q$ into the   notation.  When $x$ is chosen as in \eqref{infty x} the sum $G(x;\q)$ can be written as \eqref{infty G}, where $G_{S_\q}$ is as in \eqref{defGS}.  Then letting $F(x)$ denote the sum of $G(x;\q)$ over $\q\in\mathcal{Q}$ we obtain \eqref{low F} as before, with $\lambda_{v_0}(Y^2)$ replaced by $g_\nu(Y^2)\sim Y^{-2\beta}$.  We apply Proposition \ref{S-adic} unchanged, and arguing as in the proof of Proposition \ref{righthandside} one finds 
\begin{equation*}
F(x^2)=Y^{-2\beta}\sum_{\q\in\mathcal{Q}}L_\q(x)+O_\varepsilon(Y^{-1}+Q^{\frac{m+1}{2}+\varepsilon}).
\end{equation*}
By \eqref{newL} we obtain
\begin{equation*}
\frac{Q}{\log Q}+O_\varepsilon(Y^{2\beta}(Y^{-1-\varepsilon}+Q^{\frac{m+1}{2}+\varepsilon}))\ll\sum_{\q\in\mathcal{Q}}\frac{1}{\phi(\q)}\sum_{\chi\in X_\varepsilon (\q)} |L^S(2\beta , \pi, {\rm sym}^2 \times \chi^2) |.
\end{equation*}
We take $Y=Q^{-\frac{m+1}{2}}$.  Then for any $\beta>\frac{1}{2}-\frac{1}{m+1}$ we find
\begin{equation*}
\frac{Q}{\log Q}\ll\sum_{\q\in\mathcal{Q}}\frac{1}{\phi(\q)}\sum_{\chi\in X_\varepsilon (\q)} |L^S(2\beta , \pi, {\rm sym}^2 \times \chi^2) |.
\end{equation*}
Since the characters in $X_\varepsilon(\q)$ are all trivial at $\nu$, the standard argument of Luo-Rudnick-Sarnak then implies Theorem \ref{aux thm} for the place $\nu$.\qed

\bibliographystyle{plain}

\printindex

\end{document}